# *Chapter 8*

# Multi-Criteria Decision-Making: Aggregation-Type Methods

**Zhiyuan Wang[a], Gade Pandu Rangaiah[b, c]**

[a] School of Business, Singapore University of Social Sciences, Singapore 599494, Singapore

[b] Department of Chemical and Biomolecular Engineering, National University of Singapore, Singapore 117585, Singapore

[c] School of Chemical Engineering, Vellore Institute of Technology, Vellore 632014, India

**Abstract**

This chapter describes selected aggregation-type multi-criteria decision-making (MCDM) methods that convert an alternatives-criteria matrix (ACM) into a single performance score per alternative through additive, multiplicative or hybrid manipulations, for ranking the alternatives. The 8 methods are: Simple Additive Weighting (SAW), Multiplicative Exponent Weighting (MEW), Analytic Hierarchy Process (AHP), Analytic Network Process (ANP), Complex Proportional Assessment (COPRAS), Multi-Objective Optimization on the basis of Ratio Analysis (MOORA), Faire Un Choix Adéquat (FUCA) and Weighted Aggregated Sum Product Assessment (WASPAS). This chapter details the algorithm of each method step-by-step, illustrating every procedure with a common ACM example and full numerical calculations. Practical strengths and weaknesses of every method are outlined. A consolidated summary shows how different methods can lead to variations in the final rankings. This chapter enables the readers to: (1) explain the principles and algorithms of aggregation-type methods covered, (2) implement them on an ACM, and (3) select one or more suitable aggregation-type MCDM methods for their applications.

## 8.1 Overview

For a given alternatives-criteria matrix (ACM), multi-criteria decision-making (MCDM) serves as an effective approach for selecting the best alternative(s) from the matrix (Baydaş et al., 2024). Among the diverse MCDM methodologies, this chapter focuses on aggregation-type MCDM methods. These methods compute an overall performance score for each alternative by aggregating the values of all criteria using different principles or formulations (Wang et al., 2024). The aggregation process typically follows either an additive, multiplicative, or hybrid approach, ensuring that the final ranking of alternatives reflects the collective impact of all criteria.

This chapter covers 8 aggregation-type MCDM methods, detailing their principles and algorithms. These methods are selected based on either their widespread popularity for various applications or their recent development. Each method is elucidated through a simple numerical example, ensuring that readers gain a clear understanding of its inner workings.

The following sections in this chapter present the 8 aggregation-type MCDM methods, in chronological order. Readers need not study them sequentially and can read any of them (independent of others).

    Section 8.2    Simple Additive Weighting (SAW)
    Section 8.3    Multiplicative Exponent Weighting (MEW)
    Section 8.4    Analytic Hierarchy Process (AHP)
    Section 8.5    Analytic Network Process (ANP)
    Section 8.6    Complex Proportional Assessment (COPRAS)
    Section 8.7    Multi-Objective Optimization on the basis of Ratio Analysis (MOORA)
    Section 8.8    Faire Un Choix Adéquat (FUCA)
    Section 8.9    Weighted Aggregated Sum Product Assessment (WASPAS)
    Section 8.10   Summary

Our Microsoft Excel based program, EMCDM, described in Chapter 10, includes SAW, COPRAS, MOORA, and FUCA methods.

The learning outcomes of this chapter are: (1) Describe the principles and algorithms of aggregation-type MCDM methods covered; (2) Apply any of these methods to ACM of an application; and (3) Select one or more aggregation-type MCDM methods for solving MCDM problems.

Unless otherwise stated, the notation used throughout this chapter are as follows:
- $i \in \{1,2,\ldots,m\}$ and $j \in \{1,2,\ldots,n\}$, where $m$ is the number of rows (i.e., alternatives) and $n$ is the number of columns (i.e., criteria) in the ACM.
- $f_{ij}$ is the value of the $j^{th}$ criterion for the $i^{th}$ alternative in the original ACM.
- $F_{ij}$ is the normalized value of $f_{ij}$.
- $w_j$ is the weight assigned to the $j^{th}$ criterion, and $\sum_{j=1}^{n} w_j = 1$.
- $v_{ij}$ is the value of the $j^{th}$ criterion for the $i^{th}$ alternative in the weighted normalized ACM.

Additionally, the ACM presented in Table 8.1 is employed to walk through the steps of the aggregation-type MCDM methods covered in this chapter. This ACM is made up of 4 alternatives (i.e., $m = 4$ and $i \in \{1, 2, 3, 4\}$) and 3 criteria (i.e., $n = 3$ and $j \in \{1, 2, 3\}$). Out of the 3 criteria (C1, C2, and C3), C1 and C2 are benefit criteria to be maximized, while C3 is a cost criterion to be minimized. For illustration, the 3 criteria are assigned weights as: $w_1 = 0.25$, $w_2 = 0.33$, and $w_3 = 0.42$. All the 4 alternatives (A1, A2, A3, and A4) are non-dominated, which means that the improvement in one criterion will inevitably result in degradation in at least one other criterion (Wang et al., 2023).

Table 8.1: ACM employed for illustration of aggregation-type MCDM methods (C1 and C2 for maximization; C3 for minimization).

| Alternatives | C1 (Max) | C2 (Max) | C3 (Min) |
|---|---|---|---|
| A1 | 0.93 | 600 | 8.25 |
| A2 | 0.51 | 700 | 6.33 |
| A3 | 0.77 | 500 | 3.16 |
| A4 | 0.82 | 400 | 2.98 |
| Weight | 0.25 | 0.33 | 0.42 |

## 8.2 Simple Additive Weighting (SAW)

The SAW method, also known as the weighted sum method (WSM), described in MacCrimmon (1968), is one of the most fundamental and simple methods in MCDM. It is based on the principle that the top-ranked alternative is the one with the highest weighted sum of normalized criteria values (Nabavi et al., 2023), after converting minimization or cost criteria into maximization or benefit type. The method assumes that the criteria are compensatory, which means that a lower score in one criterion can be offset by a higher score in another. The steps in SAW are as follows.

**Step 1.** Normalize the original ACM with $m$ rows (i.e., alternatives) and $n$ columns (i.e., criteria) by using the Max normalization method. For a maximization criterion, the normalized value, $F_{ij}$ is obtained as:

$$F_{ij} = \frac{f_{ij}}{\max_{k \in \{1,2,\dots,m\}} f_{kj}} \tag{8.1}$$

For a minimization criterion, the normalized value, $F_{ij}$ is obtained as:

$$F_{ij} = \frac{\min_{k \in \{1,2,\dots,m\}} f_{kj}}{f_{ij}} \tag{8.2}$$

The above normalization ensures that all minimization criteria are converted to maximization type, indicating that a higher $F_{ij}$ always indicates a better performance (Nabavi et al., 2024).

*Numerical Calculations (using ACM from Table 8.1):*

For instance, for $i = 1$ and $j = 2$:

$$F_{12} = \frac{600}{700} = 0.8571$$

Similar calculations are performed for all other values in the ACM across the 4 alternatives and 3 criteria. The complete results, forming the normalized ACM, are shown in Table 8.2.

Table 8.2: Normalized ACM for SAW walkthrough.

| Alternatives | C1 | C2 | C3 |
|---|---|---|---|
| A1 | 1 | 0.8571 | 0.3612 |
| A2 | 0.5484 | 1 | 0.4708 |
| A3 | 0.8280 | 0.7143 | 0.9430 |
| A4 | 0.8817 | 0.5714 | 1 |

**Step 2.** Construct the weighted normalized ACM by multiplying each normalized criterion value by its corresponding weight $w_j$.

$$v_{ij} = F_{ij} \times w_j \tag{8.3}$$

*Numerical Calculations:*

For instance, for $i = 1$ and $j = 2$:

$$v_{12} = F_{12} \times w_2 = 0.8571 \times 0.33 = 0.2829$$

Likewise, calculations are carried out for all other values in the normalized ACM, using their respective assigned weights. The results, which constitute the weighted normalized ACM, are displayed in Table 8.3.

Table 8.3: Weighted normalized ACM for SAW walkthrough.

| Alternatives | C1 | C2 | C3 |
|---|---|---|---|
| A1 | 0.2500 | 0.2829 | 0.1517 |
| A2 | 0.1371 | 0.3300 | 0.1977 |
| A3 | 0.2070 | 0.2357 | 0.3961 |
| A4 | 0.2204 | 0.1886 | 0.4200 |

**Step 3.** Compute the performance score ($P_i$) for each alternative, which is obtained by summing its weighted normalized values.

$$P_i = \sum_{j=1}^{n} v_{ij} \tag{8.4}$$

Finally, the alternatives are ranked in descending order of $P_i$, with the highest value indicating the top-ranked alternative and recommended to the decision-maker.

*Numerical Calculations:*

For instance, for $i = 1$:

$$P_1 = \sum_{j=1}^{3} v_{1j} = 0.2500 + 0.2829 + 0.1517 = 0.6846$$

Similarly, the performance scores of all alternatives are calculated as $P_1 = 0.6846$, $P_2 = 0.6648$, $P_3 = 0.8388$, and $P_4 = 0.8290$. Accordingly, the ranking is A3 > A4 > A1 > A2, with A3 being the top-ranked alternative by SAW.

The advantages of SAW are as follows. (1) It is one of the simplest MCDM methods, requiring only the sum of the normalized weighted criteria values, which makes it very straightforward to understand and implement (Wang et al., 2022). (2) Its computational procedure is relatively transparent, allowing decision-makers to easily observe how each criterion and weight contributes to the overall score. (3) Due to its simplicity, it can be readily applied in many fields, such as engineering, project selection, and resource allocation, where a quick ranking of alternatives is needed.

The limitations of SAW are as follows. (1) Similar to other MCDM methods, rank reversal can occur if new alternatives are added or existing ones are removed, potentially affecting the stability of the ranking results. (2) It treats all performance values additively; high performance on one criterion is possible to completely compensate for low performance on another, which might be undesirable in certain decision-making scenarios.

## 8.3  Multiplicative Exponent Weighting (MEW)

The MEW method, also known as the weighted product method (WPM), was introduced in Miller and Starr (1969). It is an MCDM technique that evaluates alternatives based on multiplicative aggregation. Unlike the SAW method, the MEW method uses the product of normalized criteria values raised to the power of their respective weights. The alternative that has a larger product value is ranked higher.

**Step 1.** Normalize the original ACM with $m$ rows (i.e., alternatives) and $n$ columns (i.e., criteria) by using the Max normalization method. This step is the same as the Step 1 of the SAW method. For the sake of brevity, the normalization equations and numerical calculations (same as Table 8.2) are not repeated here.

**Step 2.** Calculate the performance score ($P_i$) for each alternative using the multiplicative aggregation equation as follows.

$$P_i = \prod_{j=1}^{n} F_{ij}^{w_j} \tag{8.5}$$

The alternative with the largest $P_i$ is top-ranked and recommended to the decision-maker.

*Numerical Calculations:*

For instance, for $i = 1$ and using normalized ACM in Table 8.2:

$$P_1 = \prod_{j=1}^{n} F_{1j}^{w_j} = 1.0000^{0.25} \times 0.8571^{0.33} \times 0.3612^{0.42} = 0.6197$$

Similarly, the performance scores for all alternatives are calculated as $P_1 = 0.6197$, $P_2 = 0.6271$, $P_3 = 0.8329$, and $P_4 = 0.8056$. Accordingly, the ranking is A3 > A4 > A2 > A1, with A3 being the top-ranked alternative by MEW. In this example, the top-ranked alternative, namely A3, is the same by both SAW and MEW. However, this may not always be the case in other applications.

The advantages of MEW are as follows. (1) As it uses a multiplicative aggregation, an alternative must perform reasonably well in each criterion to achieve a high overall performance score; this property helps avoid the possible excessive compensation (i.e., very high performance in one criterion fully offsets poor performance in another, as can occur in SAW). (2) Conceptually, MEW is straightforward, and decision-makers can easily understand the mechanics of how each criterion contributes to the final performance score. (3) Due to its multiplicative structure, MEW is sensitive to extreme degradation in a particular criterion, as a low value can substantially reduce the overall product. This characteristic makes it well-suited for risk-averse applications.

The limitations of MEW are as follows. (1) If the performance on a single criterion is zero (or extremely small) for a particular alternative, the product becomes zero (or very close to zero), effectively disqualifying that alternative regardless of its performance on other criteria. (2) Like other MCDM methods, rank reversal can still occur if new alternatives are introduced or existing ones are deleted, since the normalization step and multiplicative aggregation can shift the overall scale and relative scores among alternatives.

## 8.4    Analytic Hierarchy Process (AHP)

The AHP method, originally developed by Saaty (1977, 1990), is based on decision-makers' subjective pairwise comparisons. To facilitate this, Saaty's fundamental 1–9 scale (Table 8.4) is typically employed. The AHP method begins by constructing a decision hierarchy, organizing the overall goal, $n$ criteria, and $m$ alternatives into a structured model, as illustrated in Figure 8.1. Decision-makers then perform subjective pairwise comparisons to assess the relative importance among each pair of criteria using Saaty's 1–9 scale. These comparisons are compiled into an $n \times n$ pairwise comparison matrix. Following this, the eigenvalue method is used to find the priority eigenvector from the matrix, which then determines the criteria weights (Saaty, 1990). Instead of the eigenvalue method, there is a simpler approximation technique (based on sum normalization of the pairwise comparison matrix) to calculate the criteria weights, which is adopted in many applications (e.g., Teknomo, 2006; Karim & Karmaker, 2016; Wang et al., 2020). The procedure for determining criteria weights follows the same steps outlined in the AHP section of Chapter 6. If the criteria weights have already been explicitly predefined by decision-makers, this step may be omitted. Note that determination of criteria weights by the AHP method does not require ACM.

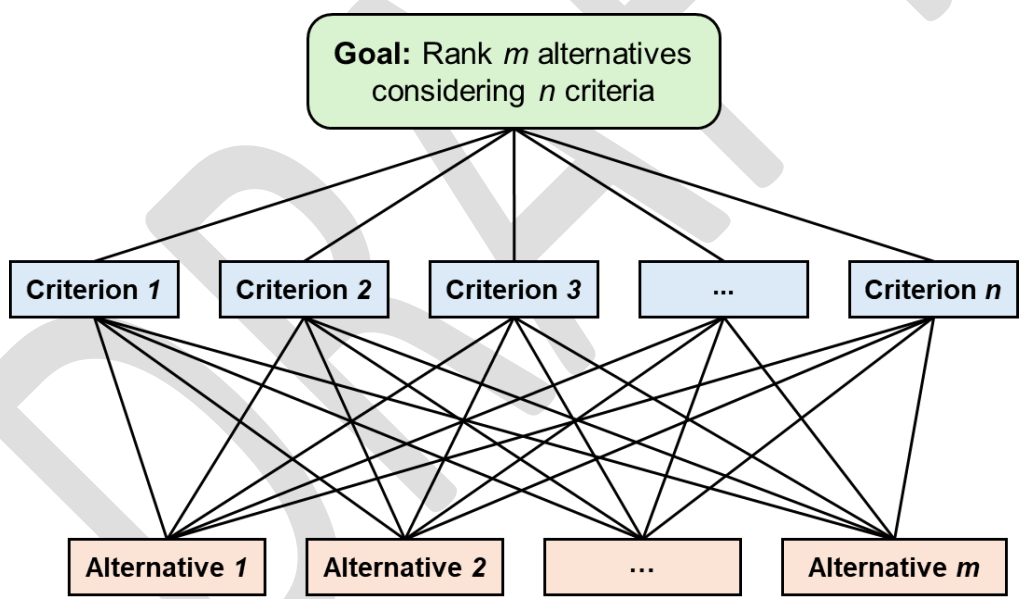

Figure 8.1: Illustration of the AHP decision hierarchy of Goal, Criteria and Alternatives in top, middle and bottom levels

Table 8.4: Saaty's 1–9 pairwise comparison scale.

| Scale | Compare between events *x* and *y* |
|:---:|:---:|
| 1 | Equal importance of *x* against *y* |
| 3 | Moderate importance of *x* against *y* |
| 5 | Essential or strong importance of *x* against *y* |
| 7 | Very strong importance of *x* against *y* |
| 9 | Extreme importance of *x* against *y* |
| 2, 4, 6, 8 | Intermediate value between adjacent scales |
| Any decimal within [1, 9], e.g., 2.3 and 5.1 | Finer scale |
| Reciprocals, e.g., 1/3 and 1/7 | If event *x* has one of the above numbers assigned to it when compared with event *y*, then *y* has the reciprocal value when compared with *x*. |

In addition to determining the criteria weights, it is crucial to recognize that the AHP method can be directly utilized for decision-making (i.e., act as an MCDM method). For this, the conventional procedure is that decision-makers conduct subjective pairwise comparisons among alternatives with respect to each criterion, leading to $n$ number of $m \times m$ pairwise comparison matrices. According to the original AHP method (Saaty, 1990), these pairwise comparisons among alternatives are purely subjective, sourced directly from the judgments of the decision-makers. For each $m \times m$ pairwise comparison matrix, AHP calculates the local priority vector (sized $m \times 1$) for alternatives with respect to the corresponding criterion using the eigenvalue method, analogous to the calculation of criteria weights at the beginning.

Notably, a typical ACM containing objective quantitative data on different measurement scales (e.g., Table 8.1), rather than the subjective pairwise comparisons using Saaty's 1–9 scale, is neither required nor immediately compatible with the AHP acting as an MCDM method. However, mapping the objective quantitative data (e.g., Table 8.1) to Saaty's 1–9 scale and constructing $n$ number of $m \times m$ pairwise comparison matrices, is feasible (Si et al., 2016).

Subsequently, AHP iterates through all $n$ criteria, calculates, and amalgamates these local priority vectors ($m \times 1$) into a local priority matrix ($m \times n$), in which each row indicates the priority of an alternative with respect to all $n$ criteria. This $m \times n$ local priority matrix is then multiplied by the $n \times 1$ criteria weights vector (previously calculated by AHP or predefined by decision-makers), resulting in an $m \times 1$ vector representing the global priorities of all alternatives. The alternative with the highest global priority value is top-ranked and recommended to decision-makers.

The main procedural steps of the AHP method are summarized below. For completeness, we will go through the full AHP steps along with numerical calculations, starting with determining criteria weights, followed by calculating the local priority vectors of alternatives with respect to each criterion, and finally aggregating the results to rank the alternatives. In the following, Steps 1 to 3 are for the criteria weights, and Steps 4 and 5 are for ranking the alternatives.

**Step 1.** Construct the $n \times n$ pairwise comparison matrix for $n$ criteria as follows.

$$\mathbf{M_A} = \begin{bmatrix} a_{11} & a_{12} & a_{13} & \cdots & a_{1n} \\ a_{21} & a_{22} & a_{23} & \cdots & a_{2n} \\ a_{31} & a_{32} & a_{33} & \cdots & a_{3n} \\ \vdots & \vdots & \vdots & \ddots & \vdots \\ a_{n1} & a_{n2} & a_{n3} & \cdots & a_{nn} \end{bmatrix} \quad (8.6)$$

- Here, $i, j \in \{1, 2, \ldots, n\}$, each entry $a_{ij}$ represents the relative importance of criterion $i$ compared to criterion $j$ by using Saaty's 1–9 scale (in Table 8.4).
- The leading diagonal entries $(a_{11}, a_{22}, a_{33}, \ldots, a_{nn})$ of this matrix are always equal to unity, that is, the relative importance of a criterion compared to itself is always 1.
- The following equation always holds true (i.e., entries below the leading diagonal are the reciprocals of the corresponding entries above it), as per the last row of Table 8.4.

$$a_{ij} = \frac{1}{a_{ji}} \quad (8.7)$$

- For instance, if we consider that the relative importance of criterion 1 versus criterion 2 is moderate (i.e., $a_{12} = 3$), then relative importance criterion 2 to criterion 1 is its reciprocal (i.e., $a_{21} = 1/3$). Thus, only $\frac{n(n-1)}{2}$ comparisons are needed to populate half of the matrix, with values of the other half obtained through simple reciprocals.

*Numerical Calculations:*

Instead of using the predefined weights presented in the last row of Table 8.1, suppose that the decision-maker provides the following subjective assessments for the 3 criteria (C1, C2, and C3) using Saaty's 1–9 scale:

- C1 is moderately more important than C2 (i.e., $a_{12} = 3$ and $a_{21} = 1/3$), meaning C1 is 3 times as important as C2, while C2 holds only 1/3 the importance of C1.
- C1 is strongly more important than C3 (i.e., $a_{13} = 5$ and $a_{31} = 1/5$).
- C2 is rated between moderate and strong in importance relative to C3 (i.e., $a_{23} = 4$ and $a_{32} = 1/4$).

Based on these assessments, the $3 \times 3$ pairwise comparison matrix is as follows:

$$\mathbf{M_A} = \begin{bmatrix} a_{11} & a_{12} & a_{13} \\ a_{21} & a_{22} & a_{23} \\ a_{31} & a_{32} & a_{33} \end{bmatrix} = \begin{bmatrix} 1 & 3 & 5 \\ 1/3 & 1 & 4 \\ 1/5 & 1/4 & 1 \end{bmatrix}$$

**Step 2.** Compute the priority vector (i.e., criteria weights in this step), using the eigenvalue method by solving the following equation.

$$\mathbf{M_A}\,\mathbf{w}^* = \lambda_{max}\,\mathbf{w}^* \tag{8.8}$$

Here, $\lambda_{max}$ is the largest eigenvalue (known as the principal eigenvalue) of the pairwise comparison matrix, and $\mathbf{w}^*$ is the corresponding principal eigenvector. The entries of this principal eigenvector are then normalized to sum to 1, yielding the criteria weights in this step, as follows.

$$w_j = \frac{w_j^*}{\sum_{k=1}^{n} w_k^*}, \quad j = 1, 2, \dots, n \tag{8.9}$$

Here, $w_j^*$ and $w_k^*$ represent the $j^{th}$ and $k^{th}$ entries of the principal eigenvector $\mathbf{w}^*$, respectively, while $w_j$ denotes the final weight of criterion $j$.

> *Numerical Calculations:*
>
> Solve the following equation using computational tools such as WolframAlpha (https://www.wolframalpha.com/input?i=eigenvalue+calculator) or Python code:
>
> $$\begin{bmatrix} 1 & 3 & 5 \\ 1/3 & 1 & 4 \\ 1/5 & 1/4 & 1 \end{bmatrix} \mathbf{w}^* = \lambda_{max} \mathbf{w}^*$$
>
> Use of WolframAlpha gives the principal eigenvalue $\lambda_{max} = 3.086$ and the corresponding principal eigenvector $\mathbf{w}^* = \begin{bmatrix} 6.6943 \\ 2.9876 \\ 1 \end{bmatrix}$. Note that when using different computational tools, the computed principal eigenvector may be scaled by a constant factor.
>
> Next, apply sum normalization to the principal eigenvector to derive the weights of the 3 criteria (C1, C2, and C3), as follows.
>
> $$w_1 = \frac{w_1^*}{\sum_{k=1}^{3} w_k^*} = \frac{6.6943}{6.6943 + 2.9876 + 1} = 0.6267$$
>
> $$w_2 = \frac{w_2^*}{\sum_{k=1}^{3} w_k^*} = \frac{2.9876}{6.6943 + 2.9876 + 1} = 0.2797$$
>
> $$w_3 = \frac{w_3^*}{\sum_{k=1}^{3} w_k^*} = \frac{1}{6.6943 + 2.9876 + 1} = 0.0936$$
>
> As expected, these weights very much depend on the values given by decision-makers, for pairwise comparison of criteria. Further, the above weights are different from the arbitrarily assumed weight values in Table 8.1.

**Step 3.** Check the consistency of the pairwise comparison matrix of $n$ criteria by calculating its consistency ratio ($CR$).

$$CI = \frac{\lambda_{max} - n}{n - 1}, \quad CR = \frac{CI}{RI} \tag{8.10}$$

Here, $CI$ is the consistency index, and $RI$ (shown in Table 8.5) is the random $CI$ value provided by Saaty (1990), who computed it by averaging the $CI$ values from numerous randomly generated $n \times n$ reciprocal matrices using Saaty's 1–9 scale. If $CR \leq 0.1$, consistency is considered acceptable. If not, decision-makers should review and revise their given values for pairwise comparison of criteria.

In addition, for a perfectly consistent pairwise comparison matrix (where each comparison $a_{ij}$ between criterion $i$ and criterion $j$ exactly equals the ratio of their resulting weights $w_i/w_j$, which also implies that the transitive property holds, i.e., $a_{ij}a_{jk} = a_{ik}$ for all criteria $i$, $j$, and $k$), Saaty (1990) proved that the largest eigenvalue of this matrix is equal to the number of criteria, i.e., $\lambda_{max} = n$.

**Note on consistency analysis**: Since we are dealing with human judgment for pairwise comparison matrix, it is inherently subjective and may sometimes be vague or inconsistent. For example, if a decision-maker prefers criterion 1 over criterion 2, and criterion 2 over criterion 3, then logically, he or she should prefer criterion 1 over criterion 3. This logical principle is known as the transitive property. However, in practice, due to cognitive biases or judgmental inconsistencies, decision-makers may unknowingly violate this property to some extent.

Table 8.5: Saaty's $RI$ values; for $n > 10$, see Podvezko (2009).

| $n$  | 1 | 2 | 3    | 4    | 5    | 6    | 7    | 8    | 9    | 10   |
|------|---|---|------|------|------|------|------|------|------|------|
| $RI$ | 0 | 0 | 0.58 | 0.90 | 1.12 | 1.24 | 1.32 | 1.41 | 1.45 | 1.49 |

*Numerical Calculations:*

In this example, with $n = 3$, refer to Table 8.5 and find the corresponding $RI = 0.58$.

$$CI = \frac{\lambda_{max} - n}{n - 1} = \frac{3.086 - 3}{3 - 1} = 0.0429$$

$$CR = \frac{CI}{RI} = \frac{0.0429}{0.58} = 0.074$$

Since the calculated $CR$ satisfies $CR \leq 0.1$, the decision-maker's subjective assessments for the 3 criteria are considered consistent and acceptable.

After obtaining the criteria weights, the following steps are for ranking the alternatives by the AHP method.

**Step 4.** Determine local priority scores of alternatives under each criterion, denoted as $v_{ij}$ with $i \in \{1,2,\ldots,m\}$ and $j \in \{1,2,\ldots,n\}$.

With respect to each of the $n$ criteria, a separate pairwise comparison matrix (size $m \times m$) is constructed, comparing all pairs of $m$ alternatives, which essentially follows the same procedure as described in Steps 1 to 3 including consistency analysis. For example, the decision-maker assesses alternative 1 against alternative 2 under the first criterion and assigns a relative importance value using Saaty's 1–9 scale. This process is repeated for all pairs of $m$ alternatives under each criterion. Once the $n$ number of $m \times m$ pairwise comparison matrices are established, the local priority vectors (representing the relative priorities of alternatives for each criterion) are computed using the eigenvalue method as in Step 2.

At this stage, if the decision-maker does not have an ACM (such as Table 8.1) based on objective measurement data, they can continue to rely entirely on their subjective judgments, as in Step 1, to construct all pairwise comparison matrices for alternatives under each criterion (e.g., A1 is moderately more important/better than A4 under C1; A3 is strongly more important/better than A2 under C3).

If an ACM based on objective measurement data is available, since it is not directly compatible with the conventional AHP process that relies on Saaty's 1–9 scale, the decision-maker needs to first map the objective data into Saaty's 1–9 scale (Si et al., 2016). This is done either through subjective assessment, leveraging the decision-maker's domain knowledge to inspect the ACM and assign pairwise comparison values of alternatives under each criterion using the 1–9 scale, or by applying mathematical mapping techniques. After this mapping, the decision-maker constructs the pairwise comparison matrices for alternatives for each criterion. One commonly used mathematical mapping technique is the logarithmic transformation, whose steps are outlined below.

Firstly, normalize the original ACM with $m$ rows (i.e., alternatives) and $n$ columns (i.e., criteria) by using the Max normalization method.

For a maximization criterion, the normalized value $F_{ij}$ is obtained as:

$$F_{ij} = \frac{f_{ij}}{\max_{k \in \{1,2,\ldots,m\}} f_{kj}} \tag{8.11}$$

For a minimization criterion, the normalized value $F_{ij}$ is obtained as:

$$F_{ij} = \frac{\min_{k \in \{1,2,\ldots,m\}} f_{kj}}{f_{ij}} \tag{8.12}$$

The above normalization transforms all minimization criteria to maximization type (i.e., a higher $F_{ij}$ always indicates a better performance).

Next, for each criterion $j$, determine the largest ratio ($r_{max,j}$) and smallest ratio ($r_{min,j}$) after comparing each alternative's normalized value with the minimum $F_{ij}$ in that criterion:

$$r_{max,j} = \frac{\max_i F_{ij}}{\min_i F_{ij}}, \quad i \in \{1,2,\ldots,m\} \text{ and } j \in \{1,2,\ldots,n\} \tag{8.13}$$

$$r_{min,j} = \frac{\min_i F_{ij}}{\min_i F_{ij}} = 1, \quad i \in \{1,2,\ldots,m\} \text{ and } j \in \{1,2,\ldots,n\} \tag{8.14}$$

Here, $r_{min,j}$ always equals 1 because it is the ratio of minimum $F_{ij}$ to itself.

Then, for each criterion $j$, the pairwise comparison between alternative $i$ and $k$ (i.e., $a_{ik}$) is mapped into Saaty's 1–9 scale by using the logarithmic transformation as follows:

$$\text{When } F_{ij} \geq F_{kj}, \quad a_{ik} = \frac{\ln\left(\frac{F_{ij}}{F_{kj}}\right) - \ln(r_{min,j})}{\ln(r_{max,j}) - \ln(r_{min,j})} \times (9 - 1) + 1, \quad i, k \in \{1,2,\ldots,m\} \tag{8.15}$$

For cases of $F_{ij} < F_{kj}$ when comparing $F_{ij}$ over $F_{kj}$, its pairwise comparison value $a_{ik}$ is obtained through simple reciprocal of the comparison of $F_{kj}$ over $F_{ij}$, since $a_{ik} = \frac{1}{a_{ki}}$.

The reason logarithmic transformation is commonly used in this context is that Saaty's 1–9 scale is meant to convey a ratio of importance (1 means equal, 3 means moderate, 9 means extreme, etc.). A logarithmic transformation preserves these multiplicative (ratio) relationships. For example, comparing $F_{ij} = 0.2$ and $F_{kj} = 0.1$, their ratio $\frac{F_{ij}}{F_{kj}}$ is 2, thus $\ln\left(\frac{F_{ij}}{F_{kj}}\right) = 0.693$; similarly, comparing $F_{ij} = 0.9$ and $F_{kj} = 0.45$, their ratio $\frac{F_{ij}}{F_{kj}}$ is also 2 and $\ln\left(\frac{F_{ij}}{F_{kj}}\right) = 0.693$. Since both comparisons yield the same logarithmic value, the pairwise comparison value $a_{ik}$ remains the same in both cases.

Once the pairwise comparison matrices for all alternatives under each criterion are constructed, the eigenvalue method equation from Step 2 can be re-used to compute the local priority scores ($v_{ij}$), which are derived by applying the sum normalization to the obtained principal eigenvector. In addition, the equation from Step 3 is re-used to calculate the consistency ratios of the matrices.

*Numerical Calculations (using ACM from Table 8.1):*

Apply Max normalization method, the original ACM is normalized and shown in Table 8.6.

Table 8.6: Normalized ACM for AHP walkthrough.

| Alternatives | C1 | C2 | C3 |
|---|---|---|---|
| A1 | 1 | 0.8571 | 0.3612 |
| A2 | 0.5484 | 1 | 0.4708 |
| A3 | 0.8280 | 0.7143 | 0.9430 |
| A4 | 0.8817 | 0.5714 | 1 |

For instance, for $j = 1$:

$$r_{max,1} = \frac{\max\limits_i F_{i1}}{\min\limits_i F_{i1}} = \frac{1}{0.5484} = 1.8235$$

$$r_{min,1} = \frac{\min\limits_i F_{i1}}{\min\limits_i F_{i1}} = 1$$

Next, under $j = 1$, let us compare A1's normalized value ($F_{11}$) over A3's normalized value ($F_{31}$), and perform logarithmic transformation:

$$a_{13} = \frac{\ln\left(\frac{F_{11}}{F_{31}}\right) - \ln(r_{min,1})}{\ln(r_{max,1}) - \ln(r_{min,1})} \times (9 - 1) + 1$$

$$= \frac{\ln\left(\frac{1}{0.8280}\right) - \ln(1)}{\ln(1.8235) - \ln(1)} \times (9 - 1) + 1$$

$$= 3.5140$$

Symmetrically, we can obtain that $a_{31} = \frac{1}{a_{13}} = \frac{1}{3.5140} = 0.2846$.

Likewise, the complete pairwise comparisons for all 4 alternatives under C1 are computed and presented in the first 5 columns of Table 8.7. The last column of Table 8.7 displays the calculated local priority scores ($v_{i1}$) for alternatives under C1, solved using the eigenvalue method described in Step 2. Additionally, the consistency ratio, $CR$ is computed as 0.0427 (within acceptable threshold), obtained using the equation in Step 3.

Table 8.7: Pairwise comparison matrix for alternatives under C1.

|  | A1 | A2 | A3 | A4 | Local priority score ($v_{i1}$) |
|---|---|---|---|---|---|
| A1 | 1 | 9.0000 | 3.5140 | 2.6762 | 0.5296 |
| A2 | 0.1111 | 1 | 0.1542 | 0.1365 | 0.0388 |
| A3 | 0.2846 | 6.4860 | 1 | 0.5441 | 0.1741 |
| A4 | 0.3737 | 7.3238 | 1.8378 | 1 | 0.2575 |

For the remaining criteria (C2 and C3), similar calculations are performed to construct the pairwise comparison matrices for alternatives under each of them and compute the corresponding local priority scores ($v_{i2}$ and $v_{i3}$) for alternatives, as presented in Table 8.8 and Table 8.9, respectively; the CR values are found to be 0.0677 and 0.0389, respectively, both within the acceptable threshold.

Table 8.8: Pairwise comparison matrix for alternatives under C2.

|  | A1 | A2 | A3 | A4 | Local priority score ($v_{i2}$) |
|---|---|---|---|---|---|
| A1 | 1 | 0.3121 | 3.6051 | 6.7958 | 0.2675 |
| A2 | 3.2042 | 1 | 5.8093 | 9.0000 | 0.5831 |
| A3 | 0.2774 | 0.1721 | 1 | 4.1907 | 0.1086 |
| A4 | 0.1471 | 0.1111 | 0.2386 | 1 | 0.0409 |

Table 8.9: Pairwise comparison matrix for alternatives under C3.

|  | A1 | A2 | A3 | A4 | Local priority score ($v_{i3}$) |
|---|---|---|---|---|---|
| A1 | 1 | 0.3245 | 0.1171 | 0.1111 | 0.0407 |
| A2 | 3.0812 | 1 | 0.1548 | 0.1445 | 0.0825 |
| A3 | 8.5392 | 6.4580 | 1 | 0.6846 | 0.3908 |
| A4 | 9.0000 | 6.9188 | 1.4608 | 1 | 0.4859 |

**Step 5.** Aggregate the local priority scores using the criteria weights calculated in Step 2, to determine the global priority scores ($P_i$) for alternatives, as follows:

$$P_i = \sum_{j=1}^{n} w_j\, v_{ij}, \quad i \in \{1,2,\ldots,m\} \tag{8.16}$$

The alternatives are ranked based on the $P_i$ values, where the alternative with the highest $P_i$ value is top-ranked.

*Numerical Calculations:*

For instance, for $i = 1$ and using the criteria weights computed in Step 2:

$$P_1 = \sum_{j=1}^{3} w_j\, v_{1j} = 0.6267 \times 0.5296 + 0.2797 \times 0.2675 + 0.09 \times 0.0936 = 0.4105$$

> Similarly, the performance scores for all alternatives are calculated as $P_1 = 0.4105$, $P_2 = 0.1951$, $P_3 = 0.1761$, $P_4 = 0.2183$. Accordingly, the ranking is A1 > A4 > A2 > A3, with A1 being the top-ranked alternative by AHP.

Alternatively, we can utilize SuperDecisions, a specialized software that has AHP and Analytic Network Process (ANP) implemented, to assist us in solving the foregoing numerical example. SuperDecisions (downloadable from https://superdecisions.com/) is a free educational software developed by the team of the creator of the methods, Professor Thomas Saaty. To begin, we construct the AHP model with the 3 hierarchical levels (namely, L1-Goal, L2-Criteria, and L3-Alternatives) within the interface of the software, as illustrated in Figure 8.2(a). Next, we input the pairwise comparison matrix for the 3 criteria (from the Numerical Calculations part of Step 1) and the pairwise comparison matrices for the 4 alternatives under each of the 3 criteria (from the first 5 columns of Table 8.7, 8.8, and 8.9) into SuperDecisions.

The software then calculates the criteria weights and the local priorities of alternatives under each criterion. These results, highlighted in the red boxes of Figure 8.2(b), match the values we previously computed manually. Readers can ignore all the other values reported in Figure 8.2(b) for now; they will understand their meanings in the subsequent ANP section. Finally, SuperDecisions determines the global priority scores of the alternatives, presented in Figure 8.2(c), which also align with our previous calculations. The original SuperDecisions file for this AHP model can be downloaded from our shared folder in Google Drive: https://drive.google.com/drive/folders/18ynIwnKaCRsptwBCrMIFVOlV9Mu_A7qM?usp=sharing

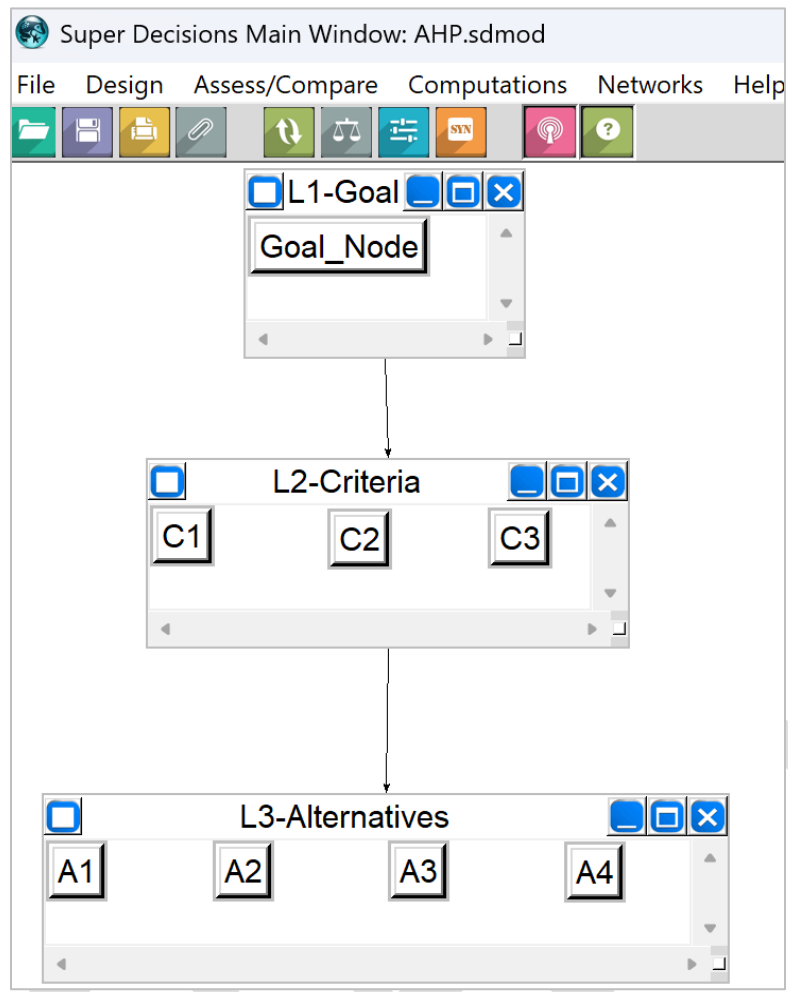

Figure 8.2(a): AHP representation in SuperDecisions software

| Cluster | Node Labels | L1-Goal | L2-Criteria | | | L3-Alternatives | | | |
|---|---|---|---|---|---|---|---|---|---|
| | | Goal_Node | C1 | C2 | C3 | A1 | A2 | A3 | A4 |
| L1-Goal | Goal_Node | 0.000000 | 0.000000 | 0.000000 | 0.000000 | 0.000000 | 0.000000 | 0.000000 | 0.000000 |
| L2-Criteria | C1 | 0.626696 | 0.000000 | 0.000000 | 0.000000 | 0.000000 | 0.000000 | 0.000000 | 0.000000 |
| | C2 | 0.279688 | 0.000000 | 0.000000 | 0.000000 | 0.000000 | 0.000000 | 0.000000 | 0.000000 |
| | C3 | 0.093616 | 0.000000 | 0.000000 | 0.000000 | 0.000000 | 0.000000 | 0.000000 | 0.000000 |
| L3-Alternatives | A1 | 0.000000 | 0.529640 | 0.267456 | 0.040746 | 0.000000 | 0.000000 | 0.000000 | 0.000000 |
| | A2 | 0.000000 | 0.038790 | 0.583045 | 0.082529 | 0.000000 | 0.000000 | 0.000000 | 0.000000 |
| | A3 | 0.000000 | 0.174103 | 0.108638 | 0.390851 | 0.000000 | 0.000000 | 0.000000 | 0.000000 |
| | A4 | 0.000000 | 0.257467 | 0.040862 | 0.485875 | 0.000000 | 0.000000 | 0.000000 | 0.000000 |

Figure 8.2(b): Criteria weights and local priorities of alternatives for each criterion inside red boxes; a screenshot from SuperDecisions for AHP.

| Here are the priorities. | | |
|---|---|---|
| Name | | Normalized by Cluster |
| A1 | | 0.41054 |
| A2 | | 0.19511 |
| A3 | | 0.17608 |
| A4 | | 0.21827 |

Figure 8.2(c): Global priorities of alternatives; a screenshot from SuperDecisions for AHP.

**Additional Comment:** If an ACM based on objective measurement data is available (e.g., Table 8.1), apart from what Step 4 describes for mapping objective data into Saaty's 1–9 scale either through subjectively assessing the ACM or mathematical transformation, another straightforward approach is to directly utilize the ACM as is and adopt a hybrid approach combining AHP with other MCDM methods (e.g., SAW, MEW) for ranking alternatives. In this hybrid approach, AHP is used only to determine criteria weights subjectively, while the objective ACM provides the evaluation data for alternatives, which are then ranked using the MCDM method (e.g., SAW, MEW). This approach, also commonly adopted in the literature (Zaidan et al., 2015; Suartini et al., 2023), bypasses the need for mapping onto Saaty's 1–9 scale and leverages objective quantitative values directly.

The advantages of AHP are as follows. (1) It provides a systematic approach to decompose complex decision problems into hierarchical levels, beginning from the overall goal, moving down to criteria, and finally to the alternatives. (2) The use of pairwise comparisons allows subjective and qualitative judgments to be included explicitly in the decision-making process while preserving a measure of logical consistency through the consistency ratio. (3) By breaking down the comparisons at each hierarchical level, AHP helps the decision-maker focus on smaller, more manageable pairs, which can reduce cognitive overload compared to simultaneously considering all criteria and alternatives. (3) Due to its intuitive nature and well-established methodology, AHP is widely adopted and/or hybridized with other methods across a broad range of decision-making contexts, including engineering design, project selection, policymaking, and resource allocation, among others.

The limitations of AHP are as follows. (1) It heavily depends on the subjective judgments of decision-makers, which can introduce biases or inconsistencies, especially when dealing with large comparison matrices. (2) The use of Saaty's 1–9 scale, while straightforward, can reduce fidelity for nuanced differences, mapping a more granular measurement scale onto these discrete intervals may oversimplify actual distinctions among criteria or alternatives. (3) Rank reversal (a common limitation of many MCDM methods) can occur if new alternatives are added or existing ones are removed, since this can alter the relative priorities derived from the pairwise comparison matrices of alternatives. (4) As the number of criteria and alternatives grows, the required pairwise comparisons increase exponentially, resulting in greater computational effort and higher cognitive demand on the part of decision-makers.

## 8.5    Analytic Network Process (ANP)

The ANP method, developed by Saaty (1996, 1999), is a generalization and an advanced version of the AHP method described in Section 8.4. This more advanced method emphasizes that everything is interconnected to everything else, with a continuous flow of influence among all elements. ANP accommodates complex decision-making scenarios involving interdependencies and feedback among criteria and alternatives. Unlike AHP, which structures decision problems strictly into hierarchical layers with unidirectional influences, ANP uses networks, allowing mutual dependence among elements and bidirectional relationships. Therefore, ANP effectively addresses decision-making complexities, such as criteria influencing alternatives and alternatives simultaneously impacting the criteria themselves.

The fundamental difference between AHP and ANP is that the latter replaces the strict hierarchical structure with a more generalized network structure comprising clusters (e.g., criteria clusters and alternatives clusters) and nodes (e.g., individual criterion and alternative). These clusters and nodes can have inner dependencies (within the same cluster) and outer dependencies (across different clusters), creating loops and feedback. The steps for applying ANP are detailed as follows:

**Step 1.** Construct the ANP network structure. For clarity of illustration, we will use a simplified example to walk through the ANP process. Suppose a company needs to select one of the 2 alternative designs ($A_1, A_2$) based on 3 criteria ($C_1$: Cost, $C_2$: Technical Feasibility, $C_3$: Market Risk). Decision-makers recognize some interdependence among these criteria and feedback from the alternatives to the criteria.

Thus, we have the following nodes:
- $\{G\}$ — single node for goal in the goal cluster.
- $\{C_1, C_2, C_3\}$ — 3 criteria nodes in the criteria cluster.
- $\{A_1, A_2\}$ — 2 alternatives nodes in the alternatives cluster.

Suppose the network of influences, as illustrated in Figure 8.3, is as follows:
- $G$ influences all criteria ($C_1, C_2, C_3$); this does not differ from that of AHP.
- $C_1$ influences $\{C_2, C_3\}$ (e.g., a higher cost can correlate with higher technical feasibility and lower market risk).
  $C_2$ influences $\{C_1, C_3\}$

  $C_3$ influences $\{C_1, C_2\}$

- Each criterion $(C_1, C_2, C_3)$ influences the alternatives $(A_1, A_2)$; this does not differ from that of AHP.
- Each alternative $(A_1, A_2)$ influences both $C_2$ and $C_3$ in a feedback loop. For example, each alternative, if implemented, might affect technical feasibility and market risk in the longer term. For instance, a more innovative design $A_1$ may push the company to adopt advanced technology solutions, raising the overall technical feasibility in the long run; and simultaneously reduce future market risk by aligning the product with emerging trends.
- If a node influences both the criteria and alternatives clusters, we assume its impact on the alternatives cluster is always 3 times as important as its impact on the criteria cluster. In other words, the particular node devotes 75% of its influence to the alternatives and 25% to the criteria.

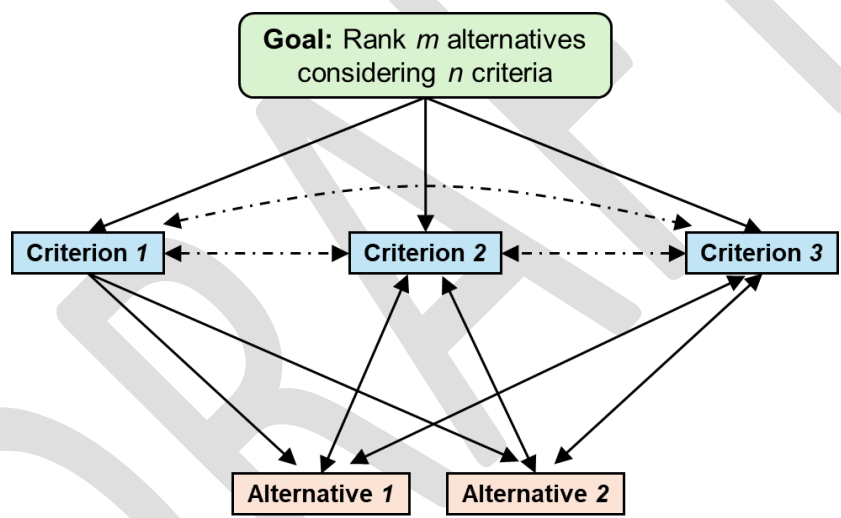

Figure 8.3: ANP representation for the simplified decision-making problem; dashed arrows indicate inner dependencies within the Criteria cluster; solid arrows represent outer dependencies across the Goal, Criteria, and Alternatives clusters. Note that some arrows are bi-directional whereas the others (e.g., from Criteria 1 to Alternatives 1 and 2) are in one direction only.

**Step 2.** Construct pairwise comparison matrices among the nodes and determine local priorities.

Using the simplified example, we now go node by node, identifying each parent node and performing the standard pairwise comparisons (as in AHP) among the child nodes under each parent, to find how strongly each parent influences its child nodes.

➤ **Under the Goal node ($G$):**

Since $G$ influences all criteria $(C_1, C_2, C_3)$:

- Parent: $G$
- Children: $\{C_1, C_2, C_3\}$

Suppose the subjective assessment using Saaty's 1–9 scale are as follows:

- $C_1$ and $C_2$ are treated equally important (ratio 1:1)
- $C_3$ is twice as important as $C_1$, and also twice as important as $C_2$ (ratios 2:1)

Hence, the $3 \times 3$ pairwise comparison matrix of the 3 criteria under $G$ is constructed as follows:

$$M_{G \to (C_1, C_2, C_3)} = \begin{pmatrix} 1 & 1 & 0.5 \\ 1 & 1 & 0.5 \\ 2 & 2 & 1.0 \end{pmatrix}$$

Same as the calculations in Step 2 of AHP, eigenvalue method is employed to determine the principal eigenvalue ($\lambda_{max}$) and normalized principal eigenvector (i.e., criteria priorities or weights). Using WolframAlpha, $\lambda_{max}$ is found to be 3.0. The corresponding criteria weights are as follows:

$$w_{G \to C_1} = 0.25, \quad w_{G \to C_2} = 0.25, \quad w_{G \to C_3} = 0.50.$$

For brevity, the consistency ratio check is omitted here. Its computation follows the exact same procedure as in Step 3 of AHP.

➢ **Under the Criterion 1 node ($C_1$):**

Since $C_1$ influences both $\{C_2, C_3\}$ and $(A_1, A_2)$:

- Parent: $C_1$
- Children: $\{C_2, C_3, A_1, A_2\}$

Besides, as aforementioned, when a node influences both the criteria and alternatives clusters, it devotes 75% of its influence to the alternatives and 25% to the criteria.

Now, inside $\{C_2, C_3\}$ pairwise comparison, $C_2$ is one-half as important as $C_3$ (ratio 1:2). The comparison matrix is as follows:

$$M_{C_1 \to (C_2, C_3)} = \begin{pmatrix} 1 & 1/2 \\ 2 & 1 \end{pmatrix}$$

For this matrix, $\lambda_{max} = 2$. The corresponding priorities are 0.3333 and 0.6667 for $C_2$ and $C_3$, respectively. Given that $C_1$ allocates 25% of its influence to criteria cluster, weights are:

$$w_{C_1 \to C_2} = 0.3333 \times 0.25 = 0.0833$$

$$w_{C_1 \to C_3} = 0.6667 \times 0.25 = 0.1667$$

Analogously, inside $\{A_1, A_2\}$ pairwise comparison, $A_1$ is 3 times as important as $A_2$ (ratio 3:1). The comparison matrix is as follows:

$$M_{C_1 \to (A_1, A_2)} = \begin{pmatrix} 1 & 3 \\ 1/3 & 1 \end{pmatrix}$$

For this matrix, $\lambda_{max} = 2$. The corresponding priorities are 0.75 and 0.25 for $A_1$ and $A_2$, respectively. Given that $C_1$ allocates 75% of its influence to the alternatives cluster, weights are:

$$w_{C_1 \to A_1} = 0.75 \times 0.75 = 0.5625$$
$$w_{C_1 \to A_2} = 0.25 \times 0.75 = 0.1875$$

➢ **Under the Criterion 2 node ($C_2$):**

The calculations under $C_2$ are similar to those under $C_1$; for brevity, only the key points and numbers are presented.

$C_2$ influences both $\{C_1, C_3\}$ and $(A_1, A_2)$:

- Parent: $C_2$
- Children: $\{C_1, C_3, A_1, A_2\}$

Inside $\{C_1, C_3\}$, $C_1$ is treated equally as important as $C_3$ (ratio 1:1). The corresponding priorities are 0.25 each. Considering 25% of parent's influence to criteria cluster, weights are:

$$w_{C_2 \to C_1} = 0.25 \times 0.25 = 0.125$$
$$w_{C_2 \to C_3} = 0.25 \times 0.25 = 0.125$$

Inside $\{A_1, A_2\}$, $A_1$ is one-half as important as $A_2$ (ratio 1:2). The corresponding priorities are 0.3333 and 0.6667 for $A_1$ and $A_2$, respectively. Considering 75% of parent's influence to the alternatives cluster, weights are:

$$w_{C_2 \to A_1} = 0.3333 \times 0.75 = 0.25$$
$$w_{C_2 \to A_2} = 0.6667 \times 0.75 = 0.50$$

➢ **Under the Criterion 3 node ($C_3$):**

Again, for conciseness, only the key points and numbers under $C_3$ are presented.

- Parent: $C_3$
- Children: $\{C_1, C_2, A_1, A_2\}$.

Inside $\{C_1, C_2\}$, $C_1$ is 3 times as important as $C_3$ (ratio 3:1). The corresponding priorities are 0.75 and 0.25 for $C_1$ and $C_2$. Considering 25% of parent's influence to criteria cluster, weights are:

$$w_{C_3 \to C_1} = 0.75 \times 0.25 = 0.1875$$

$$w_{C_3 \to C_2} = 0.25 \times 0.25 = 0.0625$$

Inside $\{A_1, A_2\}$, $A_1$ is twice as important as $A_2$ (ratio 2:1). The corresponding priorities are 0.6667 and 0.3333 for $A_1$ and $A_2$, respectively. Considering 75% of parent's influence to the alternatives cluster, weights are:

$$w_{C_3 \to A_1} = 0.6667 \times 0.75 = 0.50$$

$$w_{C_3 \to A_2} = 0.3333 \times 0.75 = 0.25$$

➤ **Under the Alternative 1 node ($A_1$):**

Since $A_1$ influences both $C_2$ and $C_3$ in a feedback loop:

- Parent: $A_1$
- Children: $\{C_2, C_3\}$.

Inside $\{C_1, C_2\}$, suppose ratio of $C_2 : C_3 = 1 : 2$, thus:

$$w_{A_1 \to C_2} = 0.3333$$

$$w_{A_1 \to C_3} = 0.6667$$

➤ **Under the Alternative 2 node ($A_2$):**

Since $A_2$ influences both $C_2$ and $C_3$ in a feedback loop:

- Parent: $A_2$
- Children: $\{C_2, C_3\}$.

Inside $\{C_1, C_2\}$, suppose ratio of $C_2 : C_3 = 1 : 1$, thus:

$$w_{A_2 \to C_2} = 0.50$$

$$w_{A_2 \to C_3} = 0.50$$

**Step 3.** Construct the weighted supermatrix (**W**).

The weighted supermatrix is formed by consolidating all the results from the previous step. It is termed "super' because it incorporates all the nodes and their influences in the ANP network. In this supermatrix, as shown in Table 8.10 each column header represents a parent node; the entries in each column indicate the influence of the parent node on its child nodes; for node pairs that were not compared in the previous step, simply assign a value of 0.

Table 8.10: ANP weighted supermatrix (**W**).

|       | $G$  | $C_1$   | $C_2$  | $C_3$   | $A_1$  | $A_2$ |
|-------|------|---------|--------|---------|--------|-------|
| $G$   | 0    | 0       | 0      | 0       | 0      | 0     |
| $C_1$ | 0.25 | 0       | 0.125  | 0.1875  | 0      | 0     |
| $C_2$ | 0.25 | 0.0833  | 0      | 0.0625  | 0.3333 | 0.50  |
| $C_3$ | 0.50 | 0.1667  | 0.125  | 0       | 0.6667 | 0.50  |
| $A_1$ | 0    | 0.5625  | 0.25   | 0.50    | 0      | 0     |
| $A_2$ | 0    | 0.1875  | 0.50   | 0.25    | 0      | 0     |

**Step 4.** Compute the limit supermatrix ($\mathbf{W}^\infty$) as follows:

$$\mathbf{W}^\infty = \lim_{k \to \infty} \mathbf{W}^k \qquad (8.17)$$

As seen, the limit supermatrix is obtained by raising the weighted supermatrix to powers until convergence. This can be computed by repeatedly multiplying the weighted supermatrix by itself until it stabilizes. In essence, the limit supermatrix represents the long-run steady-state priorities in a network of feedback influences. Each column of the supermatrix tells us how a particular parent node distributes its influence among its child nodes. By raising the supermatrix to higher and higher powers, it is essentially simulating infinitely many rounds of mutual influence among the nodes. Once the matrix powers stabilize (converge), it indicates a steady-state distribution of influence. The entries in the columns of the resulting $\mathbf{W}^\infty$, as shown in Table 8.11, can be interpreted as the overall priorities of each node under its parent node. In this example, all columns eventually converge to the same values; however, this may not be the case in all applications.

Table 8.11: ANP limit supermatrix ($\mathbf{W}^\infty$).

|       | $G$    | $C_1$  | $C_2$  | $C_3$  | $A_1$  | $A_2$  |
|-------|--------|--------|--------|--------|--------|--------|
| $G$   | 0      | 0      | 0      | 0      | 0      | 0      |
| $C_1$ | 0.0797 | 0.0797 | 0.0797 | 0.0797 | 0.0797 | 0.0797 |
| $C_2$ | 0.1991 | 0.1991 | 0.1991 | 0.1991 | 0.1991 | 0.1991 |
| $C_3$ | 0.2926 | 0.2926 | 0.2926 | 0.2926 | 0.2926 | 0.2926 |
| $A_1$ | 0.2409 | 0.2409 | 0.2409 | 0.2409 | 0.2409 | 0.2409 |
| $A_2$ | 0.1876 | 0.1876 | 0.1876 | 0.1876 | 0.1876 | 0.1876 |

Next, from $\mathbf{W}^\infty$, we extract the global priority scores for alternatives, $A_1: 0.2409$ and $A_2: 0.1876$, under the goal node $G$. After applying sum normalization, the normalized global priority score for $A_1$ is $0.5622 (= 0.2409/(0.2409 + 0.1876))$ and for $A_2$ is $0.4378$ $(= 0.1876/$ $(0.2409 + 0.1876)$. Therefore, alternative $A_1$ is ranked higher and recommended for implementation.

Alternatively, we can also use SuperDecisions to solve this problem. Figure 8.4(a) presents the weighted supermatrix ($\mathbf{W}$) calculated by SuperDecisions. Figure 8.4(b) displays the limit supermatrix ($\mathbf{W}^\infty$). Figure 8.4(c) shows the normalized global priority scores. All these results tally with our previous calculations. The original SuperDecisions file for this ANP model can be downloaded from our shared folder in Google Drive:

https://drive.google.com/drive/folders/18ynIwnKaCRsptwBCrMIFVOIV9Mu_A7qM?usp=sharing

| Cluster Node Labels | | L1-GOAL | L2-Criteria | | | L3-Alternatives | |
|---|---|---|---|---|---|---|---|
| | | Goal_Node | C1 | C2 | C3 | A1 | A2 |
| L1-GOAL | Goal_Node | 0.000000 | 0.000000 | 0.000000 | 0.000000 | 0.000000 | 0.000000 |
| L2-Criteria | C1 | 0.250000 | 0.000000 | 0.125000 | 0.187500 | 0.000000 | 0.000000 |
| | C2 | 0.250000 | 0.083333 | 0.000000 | 0.062500 | 0.333333 | 0.500000 |
| | C3 | 0.500000 | 0.166667 | 0.125000 | 0.000000 | 0.666667 | 0.500000 |
| L3-Alternatives | A1 | 0.000000 | 0.562500 | 0.250000 | 0.500000 | 0.000000 | 0.000000 |
| | A2 | 0.000000 | 0.187500 | 0.500000 | 0.250000 | 0.000000 | 0.000000 |

Figure 8.4(a): Screenshot of the weighted supermatrix ($\mathbf{W}$) from SuperDecisions for ANP.

| Cluster Node Labels | | L1-GOAL | L2-Criteria | | | L3-Alternatives | |
|---|---|---|---|---|---|---|---|
| | | Goal_Node | C1 | C2 | C3 | A1 | A2 |
| L1-GOAL | Goal_Node | 0.000000 | 0.000000 | 0.000000 | 0.000000 | 0.000000 | 0.000000 |
| L2-Criteria | C1 | 0.079748 | 0.079748 | 0.079748 | 0.079748 | 0.079748 | 0.079748 |
| | C2 | 0.199064 | 0.199064 | 0.199064 | 0.199064 | 0.199064 | 0.199064 |
| | C3 | 0.292616 | 0.292616 | 0.292616 | 0.292616 | 0.292616 | 0.292616 |
| L3-Alternatives | A1 | 0.240932 | 0.240932 | 0.240932 | 0.240932 | 0.240932 | 0.240932 |
| | A2 | 0.187639 | 0.187639 | 0.187639 | 0.187639 | 0.187639 | 0.187639 |

Figure 8.4(b): Screenshot of the limit supermatrix ($\mathbf{W}^\infty$) from SuperDecisions for ANP.

| Name | | Normalized by Cluster |
|---|---|---|
| Goal_Node | | 0.00000 |
| C1 | | 0.13956 |
| C2 | | 0.34836 |
| C3 | | 0.51208 |
| A1 | | 0.56218 |
| A2 | | 0.43782 |

Figure 8.4(c): Screenshot of the normalized global priority scores from SuperDecisions for ANP.

The advantages of ANP are as follows. (1) ANP is designed to handle interdependent criteria and alternatives by allowing feedback loops and reciprocal influences, thereby providing a more realistic representation of complex decision problems that are not strictly hierarchical. (2) The network-based model facilitates a deeper understanding of how changes to one criterion can impact others, leading to comprehensive insights that might not emerge from linear approaches. (3) By breaking down interrelationships into pairwise comparisons, ANP uses subjective judgments while preserving logical consistency to a reasonable extent. (4) Because it captures the nuanced interplay among factors, ANP can be particularly beneficial in situations where multiple variables are intertwined, and traditional methods struggle to capture those interdependencies.

The limitations of ANP are as follows. (1) ANP's feedback-driven structure can be challenging to explain to stakeholders who are unfamiliar with the method, potentially hindering consensus or agreement. (2) The method generally requires specialized software (e.g., SuperDecisions, AhpAnpLib Python library) and a solid grasp of the underlying theory, which can be a barrier in terms of both cost and the learning curve for practitioners. (3) Verification and interpretation of results can be complex, given the myriad loops and interconnections; isolating the impact of individual criteria becomes less straightforward than in hierarchical models. (4) For simpler decision problems where criteria and alternatives are relatively independent, ANP may be excessive, since more familiar and less resource-intensive methods (e.g., AHP) can often suffice while providing similar outcomes.

## 8.6    Complex Proportional Assessment (COPRAS)

The COPRAS method, proposed by Zavadskas et al. (1994), is another aggregation-type MCDM method. It employs a systematic approach to rank alternatives based on their relative importance. After constructing the weighted normalized ACM, COPRAS calculates the sum of weighted normalized values for all maximization criteria and the sum of weighted normalized values for all minimization criteria, for each alternative. These two sums are then aggregated using the method's specific equations (Zavadskas et al., 1994) to determine each alternative's performance score. The alternative with the highest performance score is top-ranked.

**Step 1.** Normalize the original ACM with m rows (i.e., alternatives) and n columns (i.e., criteria) using Sum normalization method.

$$F_{ij} = \frac{f_{ij}}{\sum_{k=1}^{m} f_{kj}} \tag{8.18}$$

*Numerical Calculations (using ACM from Table 8.1):*

For instance, for $i = 1$ and $j = 2$:

$$F_{12} = \frac{f_{12}}{\sum_{k=1}^{m} f_{k2}} = \frac{600}{600 + 700 + 500 + 400} = 0.2727$$

Similar calculations are performed for all other values in the ACM across the 4 alternatives and 3 criteria. The complete results, forming the normalized ACM, are shown in Table 8.12.

Table 8.12: Normalized ACM for COPRAS walkthrough.

| Alternatives | C1 | C2 | C3 |
| --- | --- | --- | --- |
| A1 | 0.3069 | 0.2727 | 0.3982 |
| A2 | 0.1683 | 0.3182 | 0.3055 |
| A3 | 0.2541 | 0.2273 | 0.1525 |
| A4 | 0.2706 | 0.1818 | 0.1438 |

**Step 2.** Construct the weighted normalized ACM by multiplying the values of each criterion with its assigned weight, $w_j$.

$$v_{ij} = F_{ij} \times w_j \tag{8.19}$$

*Numerical Calculations:*

For instance, for $i = 1$ and $j = 2$:

$$v_{12} = F_{12} \times w_2 = 0.2727 \times 0.33 = 0.0900$$

Likewise, calculations are carried out for all other values in the normalized ACM, using their respective assigned weights. The results, which constitute the weighted normalized ACM, are displayed in Table 8.13.

Table 8.13: Weighted normalized ACM for COPRAS walkthrough.

| Alternatives | C1 | C2 | C3 |
| --- | --- | --- | --- |
| A1 | 0.0767 | 0.0900 | 0.1672 |
| A2 | 0.0421 | 0.1050 | 0.1283 |
| A3 | 0.0635 | 0.0750 | 0.0641 |
| A4 | 0.0677 | 0.0600 | 0.0604 |

**Step 3.** For each alternative, calculate the sums of weighted normalized values for both maximization and minimization criteria, $S_{i+}$ and $S_{i-}$, separately.

$$S_{i+} = \sum_{j=1}^{g} v_{ij} \tag{8.20}$$

$$S_{i-} = \sum_{j=g+1}^{n} v_{ij} \tag{8.21}$$

Here, $g$ is the number of maximization criteria, which can be arranged in the first $g$ columns of the weighted normalized ACM. The minimization criteria are in columns $g+1$ to $n$. If all criteria are for maximization, then $g = n$; whereas $g = 0$ when all criteria are for minimization.

*Numerical Calculations:*

There are 2 maximization criteria, hence $g = 2$.
For instance, for $i = 1$:

$$S_{1+} = \sum_{j=1}^{2} v_{1j} = v_{11} + v_{12} = 0.0767 + 0.0900 = 0.1667$$

$$S_{1-} = \sum_{j=2+1}^{3} v_{1j} = v_{13} = 0.1672$$

After calculating the $S_{i+}$ and $S_{i-}$ of each alternative, the results are tabulated in the first 3 columns of Table 8.14.

Table 8.14: Weighted normalized ACM for COPRAS walkthrough.

| Alternatives | $S_{i+}$ | $S_{i-}$ | $P_i$ |
|---|---|---|---|
| A1 | 0.1667 | 0.1672 | 0.2214 |
| A2 | 0.1471 | 0.1283 | 0.2183 |
| A3 | 0.1385 | 0.0641 | 0.2813 |
| A4 | 0.1277 | 0.0604 | 0.2790 |

**Step 4.** Compute the performance score ($P_i$) of each alternative as follows.

$$P_i = \begin{cases} S_{i+} + \dfrac{\sum_{k=1}^{m} S_{k-}}{S_{i-} \sum_{k=1}^{m} \dfrac{1}{S_{k-}}}, & \text{if both maximization and minimization criteria exist} \\ S_{i+}, & \text{if only maximization criteria exist} \\ \dfrac{\sum_{k=1}^{m} S_{k-}}{S_{i-} \sum_{k=1}^{m} \dfrac{1}{S_{k-}}}, & \text{if only minimization criteria exist} \end{cases} \quad (8.22)$$

In the above equation, both $\sum_{k=1}^{m} S_{k-}$ and $\sum_{k=1}^{m} \frac{1}{S_{k-}}$ are invariant across the performance score ($P_i$) calculations for each alternative $i$, and so the term $\frac{\sum_{k=1}^{m} S_{k-}}{S_{i-} \sum_{k=1}^{m} \frac{1}{S_{k-}}}$ can be considered as a scaled $S_{i-}$ value. Further, a smaller $S_{i-}$ leads to a higher $P_i$ value. The alternative with the largest $P_i$ is top-ranked and recommended to the decision-maker.

*Numerical Calculations:*

For instance, for $i = 1$, since both maximization and minimization criteria exist:

$$P_1 = S_{1+} + \frac{\sum_{k=1}^{m} S_{k-}}{S_{1-} \sum_{k=1}^{m} \frac{1}{S_{k-}}} = 0.1667 + \frac{0.1672 + 0.1283 + 0.0641 + 0.0604}{0.1672 \left( \frac{1}{0.1672} + \frac{1}{0.1283} + \frac{1}{0.0641} + \frac{1}{0.0604} \right)} = 0.2214$$

Similarly, performance scores for all alternatives are calculated and given in the last column of Table 8.14. Accordingly, the ranking is A3 ≻ A4 ≻ A1 ≻ A2, with A3 being the top-ranked alternative by COPRAS.

The advantages of COPRAS are as follows. (1) COPRAS explicitly separates benefit (to be maximized) and cost (to be minimized) criteria, making it intuitive for decision-makers to see how each criterion contributes to or detracts from the final performance score of an alternative. (2) COPRAS has been successfully adopted in many applications (Stefano et al., 2015).

The limitations of COPRAS are as follows. (1) Because the final performance score is largely driven by the sum of weighted normalized values, exceptionally high (or low) performance on one criterion may greatly overshadow moderate performance on other criteria, which may not be desirable in some applications. (2) As with other MCDM methods, the addition or removal of alternatives may cause rank reversal.

## 8.7 Multi-Objective Optimization on the basis of Ratio Analysis (MOORA)

The MOORA method, introduced by Brauers and Zavadskas (2006), is extensively applied in various fields for MCDM. After constructing the weighted normalized ACM, the performance score for each alternative is calculated by subtracting the aggregate of minimization criteria values from the aggregate of maximization criteria values. The alternative that has a larger performance score is ranked higher. The MOORA method follows a structured sequence of the following steps.

**Step 1.** Normalize the original ACM with $m$ rows (i.e., alternatives) and $n$ columns (i.e., criteria) by Vector normalization method.

$$F_{ij} = \frac{f_{ij}}{\sqrt{\sum_{k=1}^{m} f_{kj}^2}} \tag{8.23}$$

*Numerical Calculations (using ACM from Table 8.1):*

For instance, for $i = 1$ and $j = 2$:

$$F_{12} = \frac{f_{12}}{\sqrt{\sum_{k=1}^{m} f_{k2}^2}} = \frac{600}{\sqrt{600^2 + 700^2 + 500^2 + 400^2}} = 0.5345$$

Similar calculations are performed for all other values in the ACM across the 4 alternatives and 3 criteria. The complete results, forming the normalized ACM, are shown in Table 8.15.

Table 8.15: Normalized ACM for MOORA walkthrough.

| Alternatives | C1 | C2 | C3 |
|---|---|---|---|
| A1 | 0.6015 | 0.5345 | 0.7321 |
| A2 | 0.3299 | 0.6236 | 0.5617 |
| A3 | 0.4980 | 0.4454 | 0.2804 |
| A4 | 0.5304 | 0.3563 | 0.2644 |

**Step 2.** Construct the weighted normalized ACM by multiplying the values of each criterion with its assigned weight, $w_j$.

$$v_{ij} = F_{ij} \times w_j \tag{8.24}$$

*Numerical Calculations:*

For instance, for $i = 1$ and $j = 2$:

$$v_{12} = F_{12} \times w_2 = 0.5345 \times 0.33 = 0.1764$$

Likewise, calculations are carried out for all other values in the normalized ACM, using their respective assigned weights. The results, which constitute the weighted normalized ACM, are displayed in Table 8.16.

Table 8.16: Weighted normalized ACM for MOORA walkthrough.

| Alternatives | C1 | C2 | C3 |
|---|---|---|---|
| A1 | 0.1504 | 0.1764 | 0.3075 |
| A2 | 0.0825 | 0.2058 | 0.2359 |
| A3 | 0.1245 | 0.1470 | 0.1178 |
| A4 | 0.1326 | 0.1176 | 0.1111 |

**Step 3.** Compute the performance score ($P_i$) of each alternative as follows.

$$P_i = \sum_{j=1}^{g} v_{ij} - \sum_{j=g+1}^{n} v_{ij} \tag{8.25}$$

Here, significance of $g$ is the same as that in COPRAS method, meaning that $g$ is the number of maximization criteria, which can be arranged in the first $g$ columns of the weighted normalized ACM. The minimization criteria are in columns $g+1$ to $n$. The alternative with the largest $P_i$ is top-ranked and recommended to the decision-maker. The performance score in Eq. 8.25 is somewhat simpler than that in COPRAS method (Eq. 8.22).

*Numerical Calculations:*

There are 2 maximization criteria, hence $g = 2$. For instance, for $i = 1$:

$$P_1 = \sum_{j=1}^{2} v_{1j} - \sum_{j=2+1}^{3} v_{1j} = 0.1504 + 0.1764 - 0.3075 = 0.0193$$

Similarly, the performance scores for all alternatives are calculated as $P_1 = 0.0193$, $P_2 = 0.0523$, $P_3 = 0.1537$, $P_4 = 0.1391$. Accordingly, the ranking is A3 > A4 > A2 > A1, with A3 being the top-ranked alternative by MOORA. Note that some or all $P_i$ can be negative (e.g., if all/majority of criteria are minimization type).

The advantages of MOORA are as follows. (1) It employs a straightforward mathematical formulation where the performance of each alternative is determined by summing weighted benefit criteria and then subtracting weighted cost criteria, making it easy to interpret. (2) The steps are relatively simple, which allows for efficient implementation across a wide variety of decision-making contexts.

The limitations of MOORA are as follows. (1) Similar to COPRAS, a single high or low value in one criterion can disproportionately influence the final performance score, which may be undesirable for decisions requiring a balanced performance across all criteria. (2) As with other MCDM methods, the addition or removal of alternatives may cause rank reversal.

## 8.8 Faire Un Choix Adéquat (FUCA)

FUCA is the French acronym for "Faire Un Choix Adéquat", which translates to "Make an Adequate Choice". The FUCA method, proposed by Fernando et al. (2011), determines the best alternative based on rank aggregation. It works by assigning ranks to alternatives with respect to each criterion and then computing a weighted sum of these ranks. The alternative with the smallest aggregated rank is top-ranked and recommended to the decision-maker (Wang & Rangaiah, 2017). Thus, no normalization of criteria values is required.

**Step 1.** Rank the alternatives for each criterion. Each alternative is assigned a rank with respect to every criterion, denoted as $r_{ij}$, based on its performance under that criterion, as follows:

- Rank 1 is assigned to the best alternative.
- Rank 2 is assigned to the next best alternative.
- ….
- Rank $m$ (where $m$ is the number of alternatives) is assigned to the worst alternative.

The ranking assignment depends on whether the criterion is of a benefit or cost type. For benefit criteria (maximization), the best alternative is the one with the highest value, receiving rank 1, while the worst receives rank $m$. Conversely, for cost criteria (minimization), the best alternative is the one with the lowest value, receiving rank 1, while the worst receives rank $m$.

*Numerical Calculations (using ACM from Table 8.1):*

For instance, for $i = 2$ and $j = 1$, by inspecting Table 8.1, we can easily find the rank of A2 under C1 is: $r_{21} = 4$. Similar rankings are performed for all other values in the ACM across the 4 alternatives and 3 criteria. The complete results are shown in Table 8.17.

Table 8.17: Rank of alternatives under each criterion for FUCA walkthrough.

| Alternatives | C1 | C2 | C3 |
|---|---|---|---|
| A1 | 1 | 2 | 4 |
| A2 | 4 | 1 | 3 |
| A3 | 3 | 3 | 2 |
| A4 | 2 | 4 | 1 |

**Step 2.** Compute the aggregated ranks. Once the rankings under each criterion are assigned, a weighted rank summation of all criteria ($R_i$) is performed for each alternative.

$$R_i = \sum_{j=1}^{n} (r_{ij} \times w_j) \qquad (8.26)$$

The alternative with the smallest $R_i$ is top-ranked and recommended to the decision-maker.

> *Numerical Calculations:*
>
> For instance, for $i = 4$:
>
> $$R_4 = \sum_{j=1}^{3}(r_{4j} \times w_j) = 2 \times 0.25 + 4 \times 0.33 + 1 \times 0.42 = 2.24$$
>
> Similarly, the aggregated ranks for all alternatives are calculated as $R_1 = 2.59$, $R_2 = 2.59$, $R_3 = 2.58$, $P_4 = 2.24$. Accordingly, the ranking is A4 > A3 > A1 = A2, with A4 being the top-ranked alternative by FUCA. Note that $R_i$ can have the same value for 2 or more alternatives.

The advantages of FUCA are as follows. (1) It is easy to understand and implement, requiring only ranking assignments and simple weighted summation calculations. (2) Since FUCA ranks alternatives directly instead of normalizing numerical values for each criterion, it eliminates the dependence of ranking on the normalization method. This helps avoid issues such as possibly drastic changes in normalized values when ACM is updated. (3) FUCA is less sensitive to extreme values as it considers only rankings and not absolute magnitudes.

The limitations of FUCA are as follows. (1) Since FUCA relies solely on ranks, it does not account for the actual differences between alternative performances. A small difference in criteria values may be treated the same as a large difference (i.e., in terms of ranking difference). (2) In cases where multiple alternatives have the same criterion value, additional tie-breaking rules must be defined, which can introduce subjectivity. (3) FUCA is still susceptible to rank reversal when new alternatives are introduced or existing ones are removed, as the reassignment of ranks for each criterion may alter the relative positioning of alternatives.

## 8.9    Weighted Aggregated Sum Product Assessment (WASPAS)

The WASPAS method, introduced in Zavadskas et al. (2012), is a hybrid MCDM approach that combines both the SAW and MEW methods in its algorithm. It calculates a composite performance score for each alternative by blending the additive results from SAW with the multiplicative results from MEW. By integrating these two different perspectives, WASPAS aims to leverage the strengths of each method while mitigating their individual shortcomings. The procedure of the WASPAS method consists of the following steps.

**Step 1.** Normalize the original ACM with $m$ rows (i.e., alternatives) and $n$ columns (i.e., criteria) by using Max normalization. This step is the same as the Step 1 of the SAW and MEWS methods (Sections 8.2 and 8.3). For the sake of brevity, the normalization equations and numerical calculations (same as Table 8.2) are not repeated here.

**Step 2.** Calculate the performance score ($P_i$) for each alternative by integrating additive (as in the SAW method) and multiplicative (as in the MEW method) components, as follows. The alternative with the largest $P_i$ is top-ranked and recommended to the decision-maker.

$$P_i = \lambda \sum_{j=1}^{n} F_{ij} w_j + (1-\lambda) \prod_{j=1}^{n} F_{ij}^{w_j} \qquad (8.27)$$

Here, $\lambda \in [0, 1]$ is a parameter that balances between the additive and multiplicative components. When $\lambda = 1$, WASPAS reduces to a purely additive approach (i.e., SAW). When $\lambda = 0$, it simplifies to a fully multiplicative approach (i.e., MEW). Intermediate $\lambda$ values allow a flexible combination of both approaches. Commonly, $\lambda = 0.5$ is chosen to give equal priority to both SAW and MEW components.

*Numerical Calculations:*

For instance, for $i = 1$ and $\lambda = 0.5$, using the normalized criteria values in Table 8.2:

$$P_1 = 0.5 \sum_{j=1}^{3} F_{1j} w_j + 0.5 \prod_{j=1}^{3} F_{1j}^{w_j}$$

$$= 0.5(1 \times 0.25 + 0.8571 \times 0.33 + 0.3612 \times 0.42) + 0.5(1^{0.25} \times 0.8571^{0.33} \times 0.3612^{0.42})$$

$$= 0.6521$$

Similarly, the performance scores for all alternatives are calculated as $P_1 = 0.6521$, $P_2 = 0.6460$, $P_3 = 0.8358$, $P_4 = 0.8173$. Accordingly, the ranking is A3 > A4 > A1 > A2, with A3 being the top-ranked alternative by WASPAS.

The advantages of WASPAS are as follows. (1) It involves straightforward calculations, making it easy to implement and computationally efficient. (2) By integrating SAW and MEW, WASPAS takes advantage of each method's strengths: the SAW part helps maintain clarity and proportional scaling, while the MEW part ensures no single criterion is completely overshadowed by extremely high values in others.

The limitations of WASPAS are as follows. (1) The final ranking can be affected by the choice of $\lambda$ value, which itself becomes an additional decision-making challenge. (2) If a criterion value is zero or extremely small for an alternative, the multiplicative component can excessively penalize it. (3) Similar to other MCDM methods, the addition or removal of alternatives can affect rankings, potentially leading to rank reversals.

## 8.11 Summary

Main points in this chapter on aggregation-type MCDM methods are as follows.

1. This chapter describes and illustrates 8 aggregation-type MCDM methods, each employing distinct principles for aggregating criteria values to determine an overall performance score for alternatives. Listed chronologically, these methods are SAW, MEW, AHP, ANP, COPRAS, MOORA, FUCA, and WASPAS.
2. Aggregation-type MCDM methods evaluate alternatives by combining criteria values and weights using additive, multiplicative, or hybrid aggregation techniques. These techniques transform and combine criteria values into a single overall performance score, allowing direct ranking of alternatives.
3. Each method in this chapter is systematically presented with a clear explanation of its theoretical foundation, step-by-step algorithm, and numerical calculations.
4. The chapter discusses the advantages and limitations associated with each aggregation-type MCDM method, aiding decision-makers in identifying appropriate methods aligned with specific requirements and decision contexts.
5. The hierarchical and network-based MCDM methods (i.e., AHP and ANP) offer structured approaches for incorporating subjective expert judgments through pairwise comparisons, distinguishing them from other aggregation-type methods.
6. Rank reversal, where adding or removing alternatives might alter the relative rankings of alternatives. remains a notable limitation across MCDM methods including aggregation-type.
7. A thorough understanding of the principles and relative merits of these aggregation-type MCDM methods enables decision-makers to select an appropriate method or combination of methods tailored to the specific characteristics and requirements of a decision-making scenario.
8. Finally, for methods that uses the common ACM (Table 8.1) to walk through their steps, namely, SAW, MEW, COPRAS, MOORA, FUCA, and WASPAS, it is observed that, except FUCA, all methods identify A3 as the top-ranked alternative. In contrast, FUCA ranks A4 as the top alternative. The ranking order of the remaining alternatives varies across methods, as detailed in their respective sections. These ranking results cannot be generalized, and it is better to try several MCDM methods for any application.

*Preliminary Draft Manuscript*Fernando, M. M. L., Escobedo, J. L. P., Azzaro-Pantel, C., Pibouleau, L., Domenech, S., & Aguilar-Lasserre, A. (2011, April). Selecting the best portfolio alternative from a hybrid multiobjective GA-MCDM approach for New Product Development in the pharmaceutical industry. In *2011 IEEE Symposium on Computational Intelligence in Multicriteria Decision-Making (MDCM)* (pp. 159-166). IEEE.

Karim, R., & Karmaker, C. L. (2016). Machine selection by AHP and TOPSIS methods. *American Journal of Industrial Engineering*, *4*(1), 7-13.

MacCrimmon, K. R. (1968). Decisionmaking among multiple-attribute alternatives: a survey and consolidated approach. RAND Corporation Santa Monica, https://www.rand.org/content/dam/rand/pubs/research_memoranda/2009/RM4823.pdf (accessed in March 2025).

Miller, D. W., & Starr, M. K. (1969). *Executive decision and operations research* (2nd ed.). Prentice Hall.

Nabavi, S. R., Jafari, M. J., & Wang, Z. (2023). Deep learning aided multi-objective optimization and multi-criteria decision making in thermal cracking process for olefines production. *Journal of the Taiwan Institute of Chemical Engineers*, *152*, 105179.

Nabavi, S. R., Wang, Z., & Rodríguez, M. L. (2024). Multi-Objective Optimization and Multi-Criteria Decision-Making Approach to Design a Multi-Tubular Packed-Bed Membrane Reactor in Oxidative Dehydrogenation of Ethane. *Energy & Fuels*, *39*(1), 491-503.

Podvezko, V. (2009). Application of AHP technique. *Journal of Business Economics and management*, (2), 181-189.

Saaty, T. L. (1977). A scaling method for priorities in hierarchical structures. *Journal of mathematical psychology*, *15*(3), 234-281.

Saaty, T. L. (1990). How to make a decision: the analytic hierarchy process. *European journal of operational research*, *48*(1), 9-26.

Saaty, T. L. (1996). Decisions with the analytic network process (ANP). *University of Pittsburgh (USA), ISAHP*, *96*.

Saaty, T. L. (1999). Fundamentals of the analytic network process. In *Proceedings of the 5th international symposium on the analytic hierarchy process* (Vol. 12, No. 14).

Si, J., Marjanovic-Halburd, L., Nasiri, F., & Bell, S. (2016). Assessment of building-integrated green technologies: A review and case study on applications of Multi-Criteria Decision Making (MCDM) method. *Sustainable cities and society*, *27*, 106-115.8-36

Zavadskas, E. K., Turskis, Z., Antucheviciene, J., & Zakarevicius, A. (2012). Optimization of weighted aggregated sum product assessment. *Elektronika ir elektrotechnika*, *122*(6), 3-6.

## 8.13 Exercises

E8.1 You are given an ACM consisting of 3 alternatives (A1, A2, A3) and 4 criteria (C1, C2, C3, C4), with the following values:

| Alternatives | C1 | C2 | C3 | C4 |
|---|---|---|---|---|
| A1 | 0.78 | 117 | 6.2 | 35 |
| A2 | 0.65 | 150 | 5.8 | 40 |
| A3 | 0.83 | 110 | 4.9 | 30 |

Criteria C1, C2, and C3 are benefit criteria to be maximized, whereas C4 is a cost criterion to be minimized. The weights for the criteria are: $w_1 = 0.30$, $w_2 = 0.25$, $w_3 = 0.20$, and $w_4 = 0.25$. Apply all (or some of) the aggregation-type MCDM methods covered in this chapter to rank the alternatives. Identify the top-ranked alternative by each method and compare them.

E8.2 You are given an ACM consisting of 4 alternatives (A1, A2, A3, A4) and 5 criteria (C1, C2, C3, C4, C5), with the following values:

| Alternatives | C1 | C2 | C3 | C4 | C5 |
|---|---|---|---|---|---|
| A1 | 95 | 0.245 | 13.7 | 72.1 | 0.089 |
| A2 | 88 | 0.312 | 11.9 | 69.5 | 0.075 |
| A3 | 92 | 0.298 | 12.5 | 74.8 | 0.103 |
| A4 | 85 | 0.331 | 14.1 | 70.3 | 0.091 |

Criteria C1, C2, and C4 are benefit criteria to be maximized, whereas C3 and C5 are cost criteria to be minimized. The weights for the criteria are: $w_1 = 0.20$, $w_2 = 0.20$, $w_3 = 0.25$, $w_4 = 0.15$, and $w_5 = 0.20$. Apply all (or some of) the aggregation-type MCDM methods covered in this chapter to rank the alternatives. Identify the top-ranked alternative by each method and compare them.

E8.3 You are given an ACM consisting of 6 alternatives (A1, A2, A3, A4, A5, A6) and 3 criteria (C1, C2, C3), with the following values:

| Alternatives | C1 | C2 | C3 |
|---|---|---|---|
| A1 | 5400 | 8.1 | 0.120 |
| A2 | 4750 | 7.5 | 0.135 |

| Alternatives | C1 | C2 | C3 |
|---|---|---|---|
| A3 | 6100 | 9.0 | 0.110 |
| A4 | 5900 | 8.3 | 0.095 |
| A5 | 5150 | 8.8 | 0.145 |
| A6 | 5650 | 7.9 | 0.105 |

Criteria C1 and C2 are benefit criteria to be maximized, whereas C3 is a cost criterion to be minimized. The weights for the criteria are: $w_1 = 0.40$, $w_2 = 0.35$, and $w_3 = 0.25$. Apply all (or some of) the aggregation-type MCDM methods covered in this chapter to rank the alternatives. Identify the top-ranked alternative by each method and compare them.

E8.4  Describe the similarities and differences between the SAW and MEW methods.

E8.5  Describe the similarities and differences between the AHP and ANP methods.

E8.6  Describe the similarities and differences among COPRAS, MOORA, and FUCA methods.

E8.7  Explain the role of the $\lambda$ parameter in WASPAS. Conduct sensitivity analysis by varying the value of $\lambda$ from 0 to 1 in increments of 0.1 for the ACM provided in exercise E8.1. Comment on your findings.